\documentclass{elsart}

\usepackage{amsmath,amscd,amssymb,latexsym,graphicx} 
 
\def\CC{{\mathbb{C}}}

\def\RR{{\mathbb{R}}}
\def\ZZ{{\mathbb{Z}}}

\def\cA{{\mathcal{A}}} 
 
\def\cC{{\mathcal{C}}}
\def\cE{{\mathcal{E}}}

\def\cG{{\mathcal{G}}} 
\def\cH{{\mathcal{H}}} 
\def\cK{{\mathcal{K}}} 
\def\cL{{\mathcal{L}}}
\def\cO{{\mathcal{O}}}
\def\cV{{\mathcal{V}}}
\def\cW{{\mathcal{W}}}

\DeclareMathOperator{\End}{End}   
\DeclareMathOperator{\id}{id}    

\let\le=\leqslant 
\let\ge=\geqslant \let\geq=\geqslant

\newcommand{\ot}{\otimes} 
\newcommand{\To}{\longrightarrow} 
\newcommand{\Op}{\mbox{Op}} 
\newcommand{\Mp}{\overline{M^{+}}} 
\newcommand{\co}{\overset{\text{\scriptsize{o}}}{c}}

\begin{document} 
 
\begin{frontmatter} 
 
 
\title{$K$-duality for pseudomanifolds with isolated singularities} 
 
\journal{Journal of Functional Analysis} 
 
\author{Claire Debord} 
\ead{C.Debord@math.univ-bpclermont.fr} 
\address{Laboratoire de Math{\'e}matiques, Universit{\'e} Blaise Pascal, Complexe universitaire des C{\'e}zeaux, 24 Av. des Landais, 63177 Aubi{\`e}re cedex, France} 
 
\author{Jean-Marie Lescure} 
\ead{lescure@math.univ-bpclermont.fr} 
\address{Laboratoire de Math{\'e}matiques, Universit{\'e} Blaise Pascal, 
  Complexe universitaire des C{\'e}zeaux, 24 Av. des Landais, 63177 
  Aubi{\`e}re cedex, France}

 
\begin{abstract} 
We associate to a pseudomanifold $X$ with a conical singularity a 
differentiable groupoid $G$ which plays the role of the tangent 
space of $X$. We construct a Dirac element and a Dual Dirac 
element which induce a $K$-duality between the $C^{*}$-algebras 
$C^{*}(G)$ and $C(X)$. This is a first step toward an index theory 
for pseudomanifolds. 
\end{abstract} 
 
\end{frontmatter} 
 
 
\section*{Introduction}

A basic point in the Atiyah-Singer index theory for closed 
manifolds lies in the isomorphism : 
\begin{equation}\label{symbol} 
K_*(V) \rightarrow K^*(T^*V)\ , 
\end{equation} 
induced by the map which assigns to the class of an 
elliptic pseudodifferential operator on a closed manifold $V$, 
the class of its principal symbol \cite{AS1}. 
 
To prove this isomorphism, G. Kasparov and A. Connes and G. 
Skandalis  \cite{CoS,Ka2}, define two elements $D_V\in 
KK(C(V)\ot C_0(T^*V),\CC)$ and $\lambda_V \in  KK(\CC,C(V)\ot 
C_0(T^*V))$ which induce a $K$-duality between $C(V)$ and 
$C_0(T^*V)$, i.e. : $$  \lambda _V \underset{C(V)}{\ot} D_V = 
1_{C_0(T^*V)} \ \mbox{ and } \lambda _V \underset{C_0(T^*V)}{\ot} 
D_V = 1_{C(V)} \ .$$ 
The isomorphism (\ref{symbol}) is then equal to $(\lambda 
_V \underset{C(V)}{\ot} \cdot)$. 
 
Moreover,  A. Connes and G. Skandalis recover the 
Atiyah-Singer index theorem using this $K$-duality together with 
other tools coming from bivariant $K$-theory (wrong-way functoriality maps). 
 
This notion of $K$-duality, also called Poincar\'e duality in
$K$-theory,  has a quite general meaning 
\cite{Ka1,Co0} : \\ 
two 
 $C^*$-algebras $A$ and $B$ are  $K$-dual if there exist $D\in 
 KK(A\ot B,\CC)$ and $\lambda\in KK(\CC,A\ot B)$  such that 
   $$  \lambda\underset{A}{\ot} D = 1_B\in KK(B,B) \text{ and } 
   \lambda\underset{B}{\ot} D = 1_A\in KK(A,A).$$ 
We usually call $D$ a Dirac element and $\lambda$ a dual Dirac 
element. A consequence of these equalities is that for any 
$C^*$-algebras $C$ and $E$, the groups homomorphisms : 
  $$ (\lambda\underset{A}{\ot} \cdot ): KK(A\ot C,E)\rightarrow 
  KK(C,B\ot E) $$ 
 $$ (\lambda \underset{B}{\ot} \cdot ): KK(B\ot C,E)\rightarrow 
  KK(C,A\ot E) $$ 
are isomorphisms with inverses $(\cdot\underset{B}{\ot}D)$ and 
$(\cdot\underset{A}{\ot}D)$. 
 
\smallskip It is a natural question to look for a generalization of 
the $K$-duality between a manifold and its tangent bundle for spaces less 
regular than smooth manifolds. Pseudomanifolds \cite{GoMa} offer a large class of 
interesting examples of such spaces. We have focused our attention 
on the model case of a pseudomanifold $X$ with a conical isolated 
singularity $c$. We use bivariant $K$-theory, groupoids and 
pseudodifferential calculus on groupoids to prove a Poincar\'e 
duality in $K$-theory in this context. 

Let us explain our choice of the algebras $A$ and $B$. 
As in the smooth case we take $A = C(X)$. For $B$, we need 
to define an appropriate notion of  tangent space for $X$. It 
should take into account the smooth structure of $X\setminus \{ 
c\}$ and encodes the geometry of the conical 
singularity. This problem finds an answer no longer in the 
category of vector bundles but in the larger category of groupoids. 
Thus, we assign to $X$ a smooth groupoid $G$, {\it 
the tangent space of $X$}, and we let $B$ be the non commutative 
$C^*$-algebra $C^*(G)$. 
 
The definition of the tangent space $G$ of $X$ is actually motivated
by the case of smooth manifolds. In particular, the concrete meaning
of the isomorphism (\ref{symbol}) was the initial source of inspiration.\\ The regular part 
$X\setminus \{c\}$ identifies to an open subset of $G^{(0)}$ and 
the restriction of $G$ to this subset is the ordinary tangent space of the manifold $X\setminus 
\{c\}$. The tangent space "over" the singular point is given by a 
pair groupoid. 
Furthermore the orbits space $G^{(0)}/G$ of $G$ is topologically 
equivalent to $X$, that is $C(X)\simeq C(G^{(0)}/G)$. Thus $C(X)$ 
maps to the multiplier algebra of $C^*(G)$. The Dirac element is 
defined as the Kasparov product 
 $D=[\Psi] {\ot}\partial$ where 
$[\Psi]$ is the element of $KK(C(X)\ot C^*(G), C^*(G))$ coming from 
the multiplication morphism $\Psi$ and $ \partial$ is an element of 
$KK(C^*(G),\CC)$ coming from a 
deformation groupoid $\cG$ of $G$ in a pair groupoid. This auxiliary 
groupoid $\cG$ is the analogue of the tangent 
groupoid defined by A. Connes for a smooth manifold \cite{Co0}. 
 
The construction of the dual Dirac element $\lambda$ is more 
difficult. We  let ${X_{b}}$ be  the bounded manifold with boundary 
$L$ which identifies with the closure of $X\setminus \{c \}$ in 
$G^{(0)}$, it satisfies  $X\simeq {X_{b}}/L $  and  we denote by 
$\cA G $ the Lie algebroid of the tangent space $G$. We first consider a 
suitable $K$-oriented map 
${X_{b}} \rightarrow \cA G \times {X_{b}}$. This map leads to an 
element 
$\lambda'$ of $KK(\CC, C^{*}(\cA G)\ot C(X_{b}))$. The adiabatic 
groupoid of 
 $G$ (see \cite{Co0,MP,NWX}), provides an element $\Theta$ of 
 $KK(C^{*}(\cA G),C^{*}(G))$ and we define 
$\lambda ''= \lambda' \underset{C^{*}(\cA  G)}{\ot} \Theta$. The 
element $\lambda''$ can be seen as  a continuous family 
$(\lambda''_x)_{x\in X_b}$ where $\lambda''_x\in K_0(C^*(G))$. An 
explicit description of $\lambda''$ shows  that its restriction to 
$L$ is the class of a constant family : that means $\lambda''$ 
determines an element $\lambda\in KK(\CC, C^*(G)\otimes C(X))$. 
 
An alternative and na{\"\i}ve description of $\lambda$ is the 
following. To each point $y$ of $X$ is assigned an appropriate 
open subset $\hat O_y$ of $G^{(0)}$ satisfying $K(C^*(G\vert_{\hat 
O_y}))\simeq \ZZ$. We construct a continuous family 
$(\beta_y)_{y\in X}$, where $\beta_y$ is a 
  generator of $K(C^*(G\vert_{\hat O_y}))$. 
This family gives rise to an element of $K(C^*(G\times X \vert 
_{\hat O}))$, where $\hat O$ is an open subset of $G^{(0)} \times 
X$. We obtain $\lambda$ by pushing forward this element in $K( 
C^*(G\times X))$ with the help of the inclusion morphism of 
$C^*(G\times X \vert _{\hat O})$ in $C^*(G\times X)$. 
 
The dual Dirac element has the two following important properties : 
 \begin{enumerate}\item[(i)] The set $\hat O \cap X_1 \times X_1$ is in 
   the range of the exponential map. Here $X_1$ is the complement of 
a conical open neighborhood of $c$. 
\item[(ii)] The equality $\lambda\underset{C^*(G)}{\ot} \partial = 1 
  \in K^0(X)$ holds. 
\end{enumerate} 
These two properties of $\lambda$ are crucial to obtain our main result : 
 
\smallskip \noindent {\bf{Theorem}} {\it The Dirac element $D$ and the 
  dual-Dirac element $\lambda$ induce a Poincar{\'e} duality between $C^{*}(G)$ and $C(X)$.} 
 
All our constructions are obviously equivariant under the action of a 
group of automorphisms of $X$, that is homeomorphisms of $X$ which are 
smooth diffeomorphisms of $X\setminus \{ c\}$. In previous works, 
P. Julg and G. Kasparov and G. Skandalis \cite{J,KaS} investigated the $K$-duality 
for simplicial complexes. Our approach is more in the spirit of 
\cite{Ka2} since it avoids the use of a simplicial decomposition of 
the pseudomanifold. We hope that it is better suited for applications 
to index theory. Indeed, the $K$-duality gives an isomorphism 
 between $KK(C(X),\CC)$ and 
$KK(\CC,C^*(G))$ which is, as in the smooth case, the map which 
assigns to the class of an elliptic pseudodifferential operator 
the class of its {\it symbol}. This point, among connections with 
the analysis on manifolds with boundary or conical manifolds,  will be discussed in 
a forthcoming paper. 
 
This paper is organized as follows : 
\begin{itemize} \item[ ] Section $1$ is devoted to some preliminaries 
  around  $C^*$-algebras of groupoids and special $KK$-elements. 
 \item[ ] In section $2$, we define the {\it tangent space} $G$ of a conical 
   pseudomanifold $X$ as well as the {\it tangent groupoid} $\cG$ of $X$. 
\item[ ] In section $3$, we define the Dirac element and in section 
  $4$, we construct 
the dual Dirac element. 
\item[ ] The section $5$ is devoted to the proof of the Poincar\'e 
  duality. 
\end{itemize} 
We want to address special thanks to Georges Skandalis for his always 
relevant suggestions.

\section{Preliminaries} 
 \subsection{$C^*$-algebras of a groupoid} 
We recall in this section some useful results about $C^*$-algebras of 
groupoids \cite{Re,Co0}. 
 
\smallskip \noindent Let $G \underset{r}{\overset{s}{\rightrightarrows}} G^{(0)}$ be a smooth 
Hausdorff groupoid with source $s$ and range $r$. 
If $U$ is any subset of $G^{(0)}$, we let : 
$$G_U:=s^{-1}(U)\ , \ G^U:=r^{-1}(U) \mbox{ and } 
G_U^U=G\vert_U:=G_U\cap G^U \ .$$ 
We denote by $C_c^\infty (G)$ the space of 
complex valued smooth and  compactly supported functions on $G$.  It 
is provided with a structure of involutive algebra as follows. If  $f$ and $g$ belong to 
$C_c^\infty (G)$ we define : 
 
\smallskip \noindent the {\it involution} by $$\mbox{for } \gamma\in 
G\ ,\ f^*(\gamma) = \overline{f(\gamma^{-1})} \ ;$$ 
the {\it convolution product} by $$\mbox{for } \gamma\in 
G\ ,\ f \ast 
g(\gamma)= \int_{\eta \in G^{r(\gamma)}} f(\eta)g(\eta^{-1}\gamma) \ .$$ 
 To give a sense to the integral above,  we fix a Haar 
system for $G$, that is, a smooth family $\{\lambda^x \ , \ x\in G^{(0)}\}$ of 
left invariant measures on $G$ indexed by $x\in G^{(0)}$ such that 
the support of $\lambda^x$ is $G^x$. 
 
\smallskip \noindent Alternatively, one could replace $C^\infty_c(G)$ 
by the space $C_c^\infty (G,\cL^{\frac{1}{2}})$ of compactly supported smooth sections of 
 the line  bundle of half densities $\cL^{\frac{1}{2}}$ over 
$G$. If  $k$ denotes the dimension of the $s$ (or $r$) fibers of 
$G$, the fiber $\cL^{\frac{1}{2}}_\gamma$ over $\gamma\in G$ 
is defined to be  the linear space of maps : 
$$\rho :\Lambda^k(T_\gamma(G^{r(\gamma)}))\ot\Lambda^k(T_\gamma (G_{s(\gamma)}))\rightarrow\CC$$ 
such that $\rho(\lambda v)=\vert \lambda \vert^{\frac{1}{2}} \rho(v)$ for all 
$\lambda$ in $\RR$ and $v$ in 
$\Lambda^k(T_\gamma(G^{r(\gamma)}))\ot\Lambda^k(T_\gamma(G_{s(\gamma)}))$. 
 
Then, the convolution product makes sense as the integral of a 
$1$-density on the manifold $G^{r(\gamma)}$. Both constructions lead 
to the same $C^*$-algebra.

\smallskip \noindent 
For each $x$ in $G^{(0)}$, we define a $*$-representation 
$\pi_x$ of $C_c^\infty (G)$ on the Hilbert space $L^2(G_x)$ by 
$$\pi_x(f)(\xi)(\gamma)=\int_{\eta \in G^{r(\gamma)}} 
  f(\eta)\xi(\eta^{-1}\gamma) \ ,$$ 
where $\xi \in L^2(G_x)$, $f\in C_c^\infty (G)$ and $\gamma \in G_x$.

\smallskip \noindent  The completion of $C_c^\infty (G)$ for the norm 
$\|f\|_r=\underset{x\in  G^{(0)}}{sup}\|\pi_x(f)\|$ is a 
$C^*$-algebra, called the reduced $C^*$-algebra of $G$ and denoted by 
$C^*_r(G)$.

\smallskip \noindent 
The maximal $C^*$-algebra $C^*(G)$ is the completion of 
$C_c^\infty (G)$ for the norm : 
$$\|f\|=\sup\{ \|\pi(f)\| \ | \ \pi \text{  Hilbert 
  space } *-\text{representation of } C_c^\infty (G) \} \ .$$ 
 
The previous constructions still hold when the groupoid $G$ is smooth only in the orbit direction, which means that $G\vert_{O_x}$ is smooth for any orbit $O_x=r(s^{-1}(x))$, $x\in G^{(0)}$. In this situation one can replace $C_c^\infty(G)$ by  $C_c(G)$.

\smallskip \noindent The identity map of $C_c^\infty (G)$ induces a 
surjective morphism from $C^*(G)$ onto $C_r^*(G)$. The injectivity of 
this morphism is related to amenability of groupoids 
\cite{ARe}. When $G$ is an amenable groupoid, its reduced 
and maximal $C^{*}$-algebras are equal and, moreover, this common 
$C^{*}$-algebra is nuclear. 
 
\bigskip \noindent {\it Subalgebras and exact sequences of groupoid 
  $C^*$-algebras} 
 
\smallskip \noindent 
To  an open subset  $O$  of $G^{(0)}$ corresponds an inclusion 
$i_O$ of $C_c^\infty(G\vert_O)$ into $C_c^\infty(G)$ which induces an 
injective morphism, again denoted by $i_O$, from  $C^*(G\vert_O)$ into 
$C^*(G)$. \\ When $O$ is saturated, $C^*(G\vert_O)$ 
is an ideal of $C^*(G)$. In this case, $F:=G^{(0)}\setminus O$ 
is a saturated closed subset of $G^{(0)}$ and the restriction of 
functions induces a  surjective morphism $r_F$ from $C^*(G)$ to 
$C^*(G\vert_F)$. Moreover, according to \cite{HS1}, the following sequence 
of $C^*$-algebras is exact : 
$$\begin{CD} 0 @>>> C^*(G\vert_O) @>{i_O}>> C^*(G) 
@>{r_F}>>  C^*(G\vert_F) @>>> 0 \end{CD} \ .$$ 
 
\subsection{$C^*$-modules arising from bundles and groupoids} 
Let us now consider an hermitian bundle $E$ on 
$G^{(0)}$. We equip the space $C^\infty_c(G,r^*E)$ with 
the  $C^*(G)$-valued product : 
  $$ <f , g> (\gamma) = \int_{\eta\in G^{r(\gamma)}}
  <f(\eta^{-1}),g(\eta^{-1}\gamma)>_{s(\eta)} . $$
This endows $C^\infty_c(G,r^*E)$ with a structure of 
$C^*(G)$-pre-Hilbert module and we denote by $C^*(G,E)$ the
corresponding $C^*(G)$-Hilbert module.  As usual, we note $\cL(\cE)$ and $\cK(\cE)$ the $C^*$-algebras   of (adjointable) endomorphisms and compact endomorphims of any Hilbert module $\cE$. 
\subsection{$KK$-tools} 
This paper makes an intensive use of Kasparov's bivariant 
$K$-theory. The unfamiliar reader may consult \cite{Bla,Ka1,ska}. In 
this section, we recall some basic constructions and fix the notations.

\medskip \noindent 
When $A$ is a $C^*$-algebra, the element $1_A \in KK(A,A)$ is the 
class of the triple $(A,i_A,0)$, where $A$ is graded by $A^{(1)}=0$ and 
$i_A: A \rightarrow \cL(A)$ is given by $i(a)b=ab,\ a,b\in A$.\\ 
If $B$ and $C$ are additional $C^*$-algebras, $\tau_C : 
KK(A,B)\to KK(A\ot C,B\ot C)$ is the group homomorphism defined by $\tau_C[(E,\rho,F)]=[(E\ot C,\rho \ot 
i_C , F\ot 1)]$. 
 
\smallskip \noindent The heart of Kasparov theory is the existence of 
a product which generalizes various functorial operations in
$K$-theory. Recall that the Kasparov product is a well defined
bilinear coupling $  
KK(A,B)\times KK(B,C) \rightarrow KK(A,C)$ denoted $(x,y)\mapsto 
x \underset{B}{\ot} y$ which is associative, covariant in $C$, contravariant 
  in $A$ and satisfies : 
\begin{itemize}
\item[--] $f_*(x)\underset{E}{\ot} 
  y=x\underset{B}{\ot}f^*(y)$ for any $*$-homomorphism $f:B\rightarrow E$, $x\in 
  KK(A,B)$ and $y\in KK(E,C)$.\\ 
\item[--] $x \underset{B}{\ot}1_B=1_A \underset{A}{\ot} x =x$, for $x \in KK(A,B)$.\\ 
\item[--] $\tau_D(x 
  \underset{B}{\ot} y)=\tau_D(x)\underset{B\otimes D}{\ot} \tau_D(y)$, 
  when $x\in 
  KK(A,B)$ and $y\in KK(B,C)$. 
\end{itemize} 
In the sequel, we will denote simply $x\ot y$ the product 
$x\underset{B}{\ot} y\in KK(A,C)$ when $x\in KK(A,B)$ and $y\in KK(B,C)$.\\ 
The operation $\tau_C : KK(A,B)\to KK(A\ot C,B\ot C)$ allows the 
construction of the general form of the Kasparov product : 
\begin{eqnarray} 
 KK(A_1,B_1\ot C)\times KK(A_2\ot C, B_2)\To   KK(A_1\ot A_2,B_1\ot B_2)\\ 
\label{prod} (x,y)\mapsto x\underset{C}{\ot}y:=\tau_{A_2}(x)\ot \tau_{B_1}(y) 
\end{eqnarray} 
For $x\in KK(A,B\ot C)$ and $y\in KK(B\ot C, E)$, there is an 
ambiguity in the definition of $\tau_B(x){\ot}\tau_B(y)$ : 
it can be defined by (\ref{prod}) with $B=A_2=B_1$ or by 
 $\tau_B(x\ot y)$. 
 These two products are 
different in general. Indeed,  in the 
first case, the two copies of $B$ involved in $x$ and $y$ play
different roles, 
contrary to the second case. To remove this ambiguity, we 
adopt the following convention~: 
$$\tau_B(x {\ot} y)=\tau_B(x) \ot \tau_B(y) \ \mbox{ and }$$ 
$$ x\underset{C}{\ot}y=\tau_{\underline{B}}(x)\ot \tau_B(y) \mbox{ or } 
\tau_B(x) \ot \tau_{\underline{B}}(y) \ .$$ 
Moreover, let $f_{B}:B\ot B \rightarrow B\ot B, a\ot b\mapsto b\ot a$ 
be the {\it flip} automorphism  and let $[f_{B}]$ be the corresponding element of 
$KK(B\ot B, B\ot B)$. The morphism $f_{B}$ exchanges the two copies of 
$B$, so 
$$\tau_B(x) \ot \tau_C [f_{B}]\ot \tau_B(y) = x\underset{C}{\ot}y. $$ 
 
\bigskip \noindent {\it $KK$-elements associated to deformation groupoids} 
 
\smallskip \noindent We explain here a classical construction \cite{Co0,HS1}. \\ 
A smooth groupoid $G$ is called a {\it deformation groupoid} 
if : 
$$G= G_1 \times \{0\} \cup G_2\times ]0,1] 
\rightrightarrows G^{(0)}=M\times [0,1],$$ 
where $G_1$ and $G_2$ are smooth groupoids with unit space $M$. That 
is, $G$ is obtained by gluing 
$G_2\times ]0,1]\rightrightarrows M\times ]0,1]$ 
which is the groupoid $G_2$ over $M$ parameterized by 
$]0,1]$ with the groupoid 
$G_1\times\{0\}\rightrightarrows M\times\{0\}$. 
 
\smallskip \noindent In this situation one can consider the saturated 
open subset $M\times ]0,1]$ of $G^{(0)}$. Using the isomorphisms 
$C^*(G\vert_{M\times ]0,1]}) \simeq C^*(G_2)\ot C_0(]0,1])$ and 
$C^*(G\vert_{M\times\{0\}})\simeq C^*(G_1)$,  we obtain the following 
exact sequence of $C^*$-algebras : 
$$\begin{CD} 0 @>>> C^*(G_2)\ot C_0(]0,1]) @>{i_{M\times]0,1]}}>> C^*(G) 
@>{ev_0}>>  C^*(G_1)  @>>> 0 \end{CD} 
$$ 
where $i_{M\times]0,1]}$ is the inclusion map and $ev_0$ is the {\it evaluation map} at 
$0$, that is $ev_0$ is the map coming from the 
restriction of functions to $G\vert_{M\times\{0\}}$.

We assume now that $C^*(G_1)$ is nuclear. Since 
the $C^*$-algebra $C^*(G_2)\ot C_0(]0,1])$ is contractible, 
 the long exact sequence in $KK$-theory shows that the group homomorphism 
$(ev_0)_*=\cdot {\ot}[ev_0] :KK(A,C^*(G)) \rightarrow KK(A, C^*(G_1))$ 
is an isomorphism for each $C^*$-algebra $A$. 
 
In particular with $A=C^*(G)$
we get that $[ev_0]$ is invertible in 
$KK$-theory : there is an element $[ev_0]^{-1}$ in 
$KK(C^*(G_1), C^*(G))$ such that $[ev_0] {\ot} [ev_0]^{-1}=1_{C^*(G)}$ 
and $[ev_0]^{-1} {\ot} [ev_0]=1_{C^*(G_1)}$. 
 
\smallskip \noindent Let $ev_1 :C^*(G) \rightarrow C^*(G_2)$ be the 
evaluation map at $1$ and $[ev_1]$ the corresponding element of 
$KK(C^*(G),C^*(G_2))$. 
 
\smallskip \noindent The {\it $KK$-element associated to the 
  deformation groupoid} $G$ is defined by : 
$$\delta=[ev_0]^{-1} {\ot} [ev_1]\in  KK(C^*(G_1),C^*(G_2)) \ . $$ 
\begin{exmp}\label{ex2} \begin{enumerate} 
\item  Let $G$ be a smooth groupoid and let $\cA G$ be its Lie 
  algebroid. The adiabatic groupoid of $G$ \cite{Co0,MP,NWX} : 
$$G_{ad}=\cA G\times \{0\}  \cup G\times ]0,1]   \rightrightarrows 
G^{(0)}\times [0,1], $$ 
is a deformation groupoid. Here, the vector 
bundle $\pi : \cA G\to G^{(0)}$ is considered as a groupoid in the 
obvious way.\\  Since $C_{0}(\cA^{*}G)$ 
is nuclear, the previous construction applies  and the 
 associated $KK$-element  $\delta \in KK(C_0(\cA^* G),C^{*}(G))$ gives
 rises to a map~: 
  $$  \cdot \ot \delta \ : \ K_0(C_0(\cA^*G))\To K_0(C^*(G))$$ 
This map is defined in  \cite{MP} as the analytic index of the groupoid $G$. 
 
\item  A particular case of $(1)$ is given by the tangent groupoid of 
  $\RR$ : 
$Tan(\RR)=\RR\times\RR\times]0,1]\cup 
T\RR\times\{0\}\rightrightarrows\RR\times[0,1]$ \cite{Co0}. The corresponding $KK$-element 
$\delta_B $, which belongs to $ KK(C_0(\RR^2),\cK)$ is the dual Bott element. Precisely, the map $(\cdot {\ot} \delta_B)$ induces an isomorphism from $K(C_0(\RR^2))$ into $K(\cK)\simeq \ZZ$. 
\end{enumerate}
\end{exmp} 

\subsection{Pseudodifferential calculus on groupoids}
We recall here some definitions and results of \cite{Co2,MP,NWX,Vas} (see also \cite{CoS}). 

Let $G$ be a smooth groupoid, possibly with a boundary \cite{MP}.\\ 
 Let $U_{\gamma} : C^{\infty}(G_{s(\gamma)})\To
C^{\infty}(G_{r(\gamma)})$ be the isomorphism induced by right
multiplication~:  $U_{\gamma}f(\gamma')=f(\gamma'\gamma)$. An operator
$P : C^{\infty}_{c}(G) \to C^{\infty}(G)$ is a {\it $G$-operator} if there exists a family  
$P_{x} : C^{\infty}_{c}(G_{x}) \to C^{\infty}(G_{x}) $ such that
$P(f)(\gamma)=P_{s(\gamma)}(f\vert_{G_{s(\gamma)}})(\gamma)$ and
$U_{\gamma}P_{s(\gamma)}=P_{r(\gamma)}U_{\gamma}$.  

A $G$-operator $P$ is a {\it pseudodifferential operator on $G$} (resp. of
order $m$) if for any open local chart  
$\Phi :\Omega \To s(\Omega)\times W$  of $G$  such that 
$s=pr_1\circ \Phi$ and any cut-off
function $\chi \in C^\infty_c(\Omega)$,  we have
$(\Phi^{*})^{-1}(\chi P \chi)_{x}\Phi^{*}= a(x,w,D_w)$  where 
$a\in S^*(s(\Omega)\times T^*W)$ is a classical symbol 
(resp. of order $m$).

One says that $K\subset G$ is a {\it support} of $P$ if 
$\mbox{supp}(Pf)\subset K.\mbox{supp}(f)$ for all $f\in
C^{\infty}_{c}(G)$. When $P$ has a compact support we say that $P$ is
{\it uniformly supported}.   

These definitions extend immediately to the case of  operators acting
on sections of bundles on $G^{(0)}$ (pulled back to $G$ with $r$), and
we denote by $\Psi^*(G,E)$ the algebra of uniformly supported
pseudodifferential operators on $G$ acting on sections of $E$. Thanks
to the invariance property, each operator $P\in\Psi^*(G;E) $ has a
principal symbol  
$\sigma(P)\in C^{\infty}_{c}(S^{*}G,\hom \pi^{*}E)$ where $S^*G$
is the sphere bundle associated to  $\cA^*G$  and $\pi$ its natural
projection onto $G^{(0)}$.  
The following inclusions hold : $\Psi^{0}(G,E)\subset \cL(C^*(G,E))$ and 
$\Psi ^{-1}(G,E)\subset \cK(C^*(G,E))$. Moreover, the symbol map
extends by continuity and gives rise to the following exact sequence of
$C^*$-algebras~:  
\begin{equation}\label{E-S-PDO}
 0\To \cK(C^*(G,E))\To \Psi_0(G,E)\stackrel{\sigma}{\To} C_0(S^*G,\hom \pi^*E) \To 0. 
\end{equation}
where $\Psi_0(G,E) =\overline{\Psi^0(G,E)}^{ \cL(C^*(G,E))}$.
Finally a linear section $\Op_G $ of the symbol map can be defined by the following formula~:
$$\mbox{Op}_G(a)(u)(\gamma) =\int_{ \xi\in 
    \cA^*_{r(\gamma)}G;\atop \gamma'\in G_{s(\gamma)}} 
 e^{i<E_G^{-1}(\gamma'\gamma^{-1}),\xi>}a(r(\gamma),\xi) 
 \phi(\gamma'\gamma^{-1})p_{r(\gamma'),r(\gamma)}u(\gamma')d\gamma'd\xi.$$ 
Here $\phi\in  C^\infty(G)$ is supported in the range of an
exponential map $E_G :\cV(\cA G)\to G$ where
$\cV(\cA G)$ denotes a small neighborhood of the zero section in
$\cA G$; moreover $\phi$ is assumed to be equal to
one on a neighborhood of $V$  in $G$.  We have used a parallel  transport $p$ to 
get local trivializations of the bundle $E$. It is implicit in the
formula that the symbol $a\in C^{\infty}_{c}(S^* G,\hom \pi^*E)$ has
been extended in the usual way to $\cA^{*} G$. 

\section{The Geometry} 
Let $X_1$ be an $m$-dimensional compact manifold with 
boundary $L$. We attach to each connected component $L_i$ of the 
boundary the cone  $cL_i= L_i\times [0,1] / L_i\times \{0\}$, using the 
obvious  map  $L_i\times\{1\} \to L_i\subset\partial X_1$. 
The new space $X=\sqcup_{i=1}^{p} cL_i \cup X_1$ is a compact pseudomanifold 
with isolated singularities \cite{GoMa}. In general, there is no 
manifold structure around the vertices of the cones. From now on, we 
assume  that $L$ is connected, i.e. $X$ has only one singularity 
denoted by $c$. The general case follows by exactly the same methods. 
 
 For any $\varepsilon \in ]0,1]$, we will refer to $c_\varepsilon
 L=L\times [0,\varepsilon] / L\times \{0\}$ as a compact cone over
 $L$, to 
 $\co_\varepsilon L=L\times [0,\varepsilon[ / L\times \{0\}$ as an
 open cone over $L$ and we let 
 $X_\varepsilon=L\times[\varepsilon,1]\cup X_1$.

\smallskip \noindent 
We define the manifold $M$ by attaching to $X_1$ a cylinder 
$L\times ]-1,1]$. We fix on $M$ a Riemannian metric which is of 
product type on $L\times ]-1,1]$ and we assume that its injectivity 
radius is bigger than $1$.\\ 
  We will use the following notations :  $ M^+$ denotes 
$L\times ]0,1] \cup X_1$, $\overline{M_+}$ its closure in $M$ and 
$ M^-=L\times ]-1,0 [$.  If $y$ is a point of the cylindrical part of $M$ or $X\setminus\{c\}$, we 
will write $y=(y_L,k_y)$ where $y_L \in L$ and 
$k_y\in ]-1,1] $ are the tangential and  radial 
coordinates. We extend the map $k$ on $X_1$ to a smooth defining function for its 
boundary; in particular, $k^{-1}(1)=\partial X_1$ and $k(X_1) \subset[1,+\infty[$. 
\begin{center} 
\includegraphics[width=11cm]{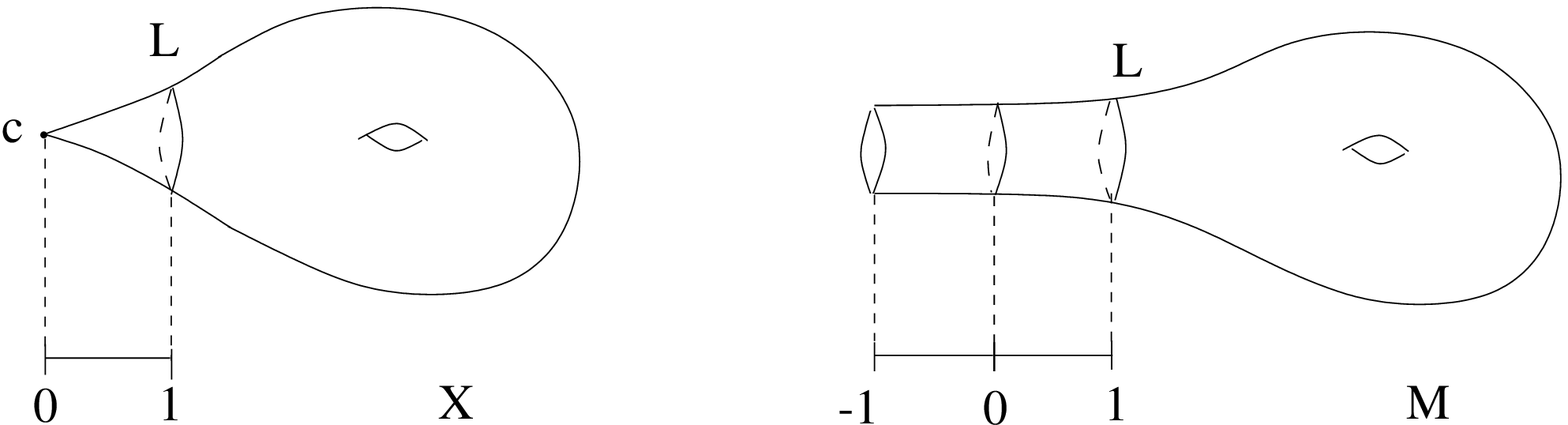} 
\end{center} 
 \subsection{The tangent space of the conical pseudomanifold $X$} 
Let us consider $T\overline{M^+}$, the restriction to $\overline{M^+}$ 
of the tangent bundle of $M$. As a $\cC^\infty$ vector bundle, it is a 
smooth groupoid with unit space $\overline{M^+}$. We define the 
groupoid $G$ as the disjoint union : 
$$G=M^-\times M^-\ \cup \ 
T\overline{M^+} \overset{s}{\underset{r}{\rightrightarrows}} M,$$ 
where  $M^-\times M^-\rightrightarrows M^-$ is the pair groupoid. 
 
\smallskip \noindent In order to endow $G$ with a smooth structure, 
compatible with the usual smooth structure on $M^-\times M^-$ and 
on $TM^+$, we have to take care of what happens around points of 
$T\overline{M^+}|_{\partial\Mp}$. \\ Let $\tau$ be a 
smooth positive function on $]-1,+\infty[$ such that
$\tau^{-1}(\{0\})=\RR^{+}$. 
We let $\tilde{\tau}$ be the smooth 
map from $M$ to $\RR^{+}$ given by $\tilde{\tau}(y)=\tau(k_y)$. \\ 
Let $(U,\phi)$ be a local chart for 
$M$ around $z\in\partial\Mp$. Setting $U^-=U\cap M^-$ and 
$\overline{U^+}= U \cap \Mp$, we define a local chart of $G$ by : 
$$ \begin{array}{cccc} \tilde{\phi} : & 
 U^-\times U^- \cup \ T \overline{U^+} & \longrightarrow &  \RR^{m}\times \RR^{m} \end{array}$$ 
\begin{equation}\label{local-charts} \tilde{\phi}(x,y)= 
(\phi(x),\frac{\phi(y)-\phi(x)}{\tilde{\tau}(x)}) \mbox{ if } (x,y)\in 
U^-\times U^- \mbox{ and} 
\end{equation} 
$$\tilde{\phi}(x,V)=(\phi(x),(\phi)_*(x,V)) \mbox{ elsewhere.}$$ 
Let us explain why the range of $\widetilde{\phi}$ is open. We can
assume that $\phi(U)=\RR^m$ and
$\phi(U^-)=\RR^m_- = \RR^{m-1}\times]-\infty,0[$.   
Let $B$ be a open ball in $\RR^m$. Since $\widetilde{\tau}$
vanishes on $\overline{U^+}$  there exists an open neighborhood $W$ of
$\partial\overline{U^+}$ such that 
$\{ \widetilde{\tau}(x)p+\phi(x) \ | \ x\in W\cap U^-,\ p\in
B\}\subset\RR^m_-$.  Then
$\widetilde{\phi}(T\partial\overline{U^+})\cap \RR^m\times B \subset
\phi(W)\times B \subset\mbox{Im}\widetilde{\phi}$.  

We define in this way a structure of smooth  groupoid on $G$. 
\begin{rem}\begin{enumerate} \item If $\tau$ is $C^l$  then the atlas
    defined above provides $G$ with a structure of $C^l$
    groupoid (it is easy to see that the source, target and inversion
    maps have the same regularity as the atlas). 
\item At the topological level, the space of orbits $M/G$ of $G$ is 
equivalent to $X$ : there is a canonical isomorphism 
between the algebras $C(X)$ and $C(M/G)$. 
\end{enumerate} 
\end{rem} 
 \begin{defn}\label{esptangent}  The smooth groupoid 
  $G \rightrightarrows M$ is called a tangent space of $X$. 
\end{defn} 
 It is important to remark that the Lie algebroid of 
$G \rightrightarrows M$ is the bundle 
$\cA G =TM$ over $M$ with anchor $p_{G}:\cA G=TM \rightarrow TM \ , \ 
(x,V) \mapsto (x,\tilde{\tau}(x)V)$; in particular $p_{G}$ is the 
zero map in restriction to $T\overline{M^+}$. The exponential 
map $exp$ of the Riemannian manifold $M$ provides an exponential map $E_G$ 
for the groupoid $G$ (for a description of 
exponential maps for groupoids, see e.g. \cite{NWX,CF}). More
precisely :  
$$ E_G : \mathcal{V}(TM) \longrightarrow  G $$ 
$$ E_{G}(y,V)=(y,V) \mbox{ when } y\in \overline{M^+} \mbox{ and }$$ 
$$E_{G}(y,V)=(y,exp_{y}(-\tilde{\tau}(y)V)) \mbox{ when } 
y\in M^-$$ where $\mathcal{V}(TM)=\{(y,V)\in TM\ | \ \|\tilde{\tau}(y)V\|<1 \mbox{ and } 
 exp_{y}(-\tilde{\tau}(y)V))\in M^-  \mbox{ if } 
y\in M^-\}.$ 
 
The map $E_G$ is a diffeomorphism onto a neighborhood of the 
unit space $M$ in $G$. In fact, we could have defined the smooth structure of $G$ using the 
map $E_{G}$. 
\begin{rem}\begin{enumerate} \item  There exists a slightly different groupoid which 
could naturally play the role of the tangent space of $X$. We will call 
it the tangent space with tail. It is defined by 
: $$G_q= L\times L \times T(]-1,0[) \cup T\overline{M^+} 
\rightrightarrows M \ .$$ 
As a groupoid, $G_q$ is the union of two groupoids : the bundle 
$T\Mp\rightrightarrows \Mp$ and the groupoid $L\times L \times 
T(]-1,0[) \rightrightarrows L\times ]-1,0[=M^-$ which is the 
product of the pair groupoid over $L$ with the vector bundle 
$T(]-1,0[)\rightrightarrows ]-1,0[$. One can equip $G_q$ with a smooth 
structure similarly as we did for $G$.
We will see that the $C^*$-algebras of $G$ and $G_q$ are
$KK$-equivalent . 
\item The groupoid $G$ is obtained by gluing along their common boundary
  $TL\times \RR$ the groupoids $T\Mp$ and a groupoid isomorphic to
  $\mbox{Tan}(L)\rtimes \RR^*_+$
  obtained by the action of $\RR_+^*$ (by multiplication on the real
  parameter) on the tangent groupoid
  $\mbox{Tan}(L)=L\times L \times \RR^*_+\cup TL\times \{0\}$ of $L$. The
  groupoid $G_q$ is defined in the same way except that we consider
  the trivial action of $\RR^*_+$. 
\end{enumerate} \end{rem} 
 \subsection{The tangent groupoid of the pseudomanifold $X$} 
The following construction is a natural generalization of the tangent groupoid of a manifold defined 
by A. Connes \cite{Co0}. We define the {\it tangent groupoid} $\cG$ of the 
pseudomanifold $X$ as a deformation of the pair groupoid 
over $M$ into the groupoid $G$. This deformation process has a nice 
description at the level of Lie algebroids. Indeed, the  Lie algebroid 
of $\cG$ should be the (unique) Lie algebroid given 
by the fiber bundle $\cA \mathcal{G}=[ 0,1] \times \cA G=[ 0,1] \times TM$ over $[ 0,1] \times M$, 
with anchor map 
$$\begin{array}{cccc} p_{\mathcal{G}}: & \cA\mathcal{G}=[ 0,1] \times TM & \longrightarrow & T([ 0,1] 
\times M)=T[ 0,1]\times TM \\ & (\lambda,x,V) & \mapsto & 
(\lambda,0,x,p_{G}(x,V)+\lambda V)=(\lambda,0,x,(\tilde{\tau 
}(x)+\lambda)V)  \ . \end{array}$$ 
Such a Lie algebroid is  almost injective, thus it is integrable 
\cite{CF,Moa3}.\\ We now define the {\it tangent groupoid} : 
$$\mathcal{G}=M\times M\times ]0,1] \ \cup \ G\times \{0\} 
{\rightrightarrows} M \times 
[0,1], $$ 
whose smooth structure is described hereafter. \\ 
Since  $\cG=(M\times M\times [0,1]\setminus(\Mp\times M \cup M \times \Mp)\times \{0\})\cup T\Mp\times\{0\}$, 
we keep the smooth structure on 
$M\times M\times[0,1]\setminus (\overline{M^+}\times M \cup M \times \overline{M^+}) \times \{0\} $ 
 as an open subset in the manifold 
with boundary $M\times M\times[0,1]$. We 
consider the following map : 
$$ \rho : \mathcal{V}(TM\times [0,1]) 
\longrightarrow \mathcal{G}=M\times M\times ]0,1] \ \cup \ G\times 
\{0\} $$ 
$$ \rho(z,V,\lambda)= \left\{ \begin{array}{lll} (z,V,0) \mbox{ if } z\in 
\overline{M^+} \mbox{ and } \lambda=0 \\ 
(z,exp_{z}(-(\tilde{\tau}(z)+\lambda).V),\lambda) \mbox{ elsewhere ,} 
\end{array} \right. $$ 
where $\mathcal{V}(TM\times [0,1]) $ is an 
open subset in $TM\times [0,1]$ such that $\mathcal{V}(TM\times [0,1]) \cap
TM\times \{0\}=\mathcal{V}(TM)$, and 
which is small enough so the exponential in the definition of $\rho$ 
is well defined. Then $T\overline{M^+}\times\{0\}$ is in the image of
$\rho$ and we equip $\cG$  around $T\overline{M^+}\times\{0\}$ with
the smooth structure for which $\rho$ is a diffeomorphism onto
its image. One can easily check that it is compatible with the smooth
structure of $M\times M\times[0,1]\setminus (\overline{M^+}\times M
\cup M \times \overline{M^+}) \times \{0\} $.   

 \noindent The Lie algebroid of $\mathcal{G}$ is $\cA 
\mathcal{G}$ and $\rho$ is an exponential map for $\mathcal{G}$. 
\subsection{ The $C^*$-algebras} 
Let $\mu$ be the Riemannian measure on $M$  and let $\nu$ be 
the corresponding Lebesgue measure on the fibers of $TM$. The family $\{\lambda^x \ ; \ x\in M\}$, where 
$\displaystyle d\lambda^x(y)=\frac{1}{\tilde{\tau}(y)^{m}} d\mu (y)$ if $x$ 
belongs to $M^-$, and $d\lambda^x(V)=d\nu (V)$ if $x$ is in 
$\overline{M^+}$, is a  Haar system for $G$. We use this Haar system 
to define the convolution algebra of $G$. 
\begin{rem}\begin{enumerate} \item \label{min=max}  
$G$ is a continuous field of amenable groupoids parameterized by $X$. More
precisely, $G=\sqcup_{x\in X} \pi^{-1}(x)$ where $\pi : G\to X$ is the
obvious projection map.  If $x\not= c$, $\pi^{-1}(x)=T_xM$ is
amenable. If $x=c$,  $\pi^{-1}(c)= M^-\times M^-\cup T\Mp\vert_{\partial \Mp}$ is
isomorphic to the groupoid $H= \mbox{Tan}(L)\rtimes \RR$ of an action
of $\RR$ on the tangent groupoid $\mbox{Tan}(L)= L\times L\times]0,1[\cup
TL\times\{0\}$ of $L$. The groupoid $H$ is an extension of the group
$\RR$ by $\mbox{Tan}(L)$, both of them being amenable, according to
\cite{ARe}, theorem 5.3.14, $H$ is amenable. 
Finally according to  \cite{ARe}, proposition 5.3.4, $G$ is  amenable. 
In the same way, $\cG$ and $G_q$ are 
amenable. Hence their reduced and  maximal $C^*$-algebras  are 
  equal and they are nuclear. 
\item Using the $KK$-equivalence between $\cK$ and  $C_0(\RR^2)$ (cf. example \ref{ex2} (2)) one
  can establish a $KK$-equivalence between $C^*(G)$ and $C^*(G_q)$. 
\end{enumerate} 
\end{rem} 
\section{The Dirac element} 
The tangent groupoid $\cG \rightrightarrows M\times[0,1]$ 
is a deformation groupoid and its $C^{*}$-algebra is nuclear, thus it 
defines a $KK$-element. 
We let $\widetilde{\partial}$ be the $KK$-element associated to 
$\cG$. More precisely : 
$$\widetilde{\partial} =[e_0]^{-1} \otimes [e_1] \in KK(C^*(G),\cK ), $$ 
where $e_0: C^*(\cG) \to C^*(\cG\vert_{M\times\{0\}})=C^*(G)$, the
evaluation map at $0$ is $K$-invertible, and 
$e_1:C^*(\cG)\to C^*(\cG\vert_{M\times\{1\}})=\cK(L^2(M))$ is the
evaluation map at $1$. Let $b$ be the (positive) generator of $KK(\cK,\CC)\simeq \ZZ$. 
We set $\partial = \widetilde{\partial}\ot b$. 
 
 \noindent The algebra $C(X)$ is isomorphic to the algebra 
of continuous functions on the orbits space $M/G$ of $G$. Thus $C(X)$ 
maps to the 
multiplier algebra of $C^*(G)$ and we let $\Psi$ 
be the morphism 
$\Psi :C^*(G) \otimes C(X) \rightarrow C^*(G)$ induced by the 
multiplication. In other words, if $a\in C_c^\infty(G)$ and $f\in C(X)$, 
$\Psi(a,f)\in C_c(G)$ is defined by 
$$\Psi(a,f)(\gamma)=\left \{ 
  \begin{array}{ll} a(\gamma)f(r(\gamma))= a(\gamma)f(s(\gamma)) \mbox{ if } \gamma\in T\Mp 
    \\ a(\gamma)f(c) \mbox{ if } \gamma\in M^-\times M^- \end{array} 
\right.$$ 
We denote by $[\Psi]$ the corresponding element in 
$KK(C^*(G)\otimes C(X),C^*(G))$. 
\begin{defn} The Dirac element is : 
$$D=[\Psi]\otimes \partial \in KK(C^*(G)\otimes C(X),\CC) \ .$$ 
\end{defn} 
\section{The dual Dirac element}  
We first recall the construction of the  dual Dirac element for a
compact manifold $V$
\cite{CoS,Ka2}.\\  
Let $V$ be a smooth compact $n$ dimensional Riemannian manifold, whose 
injectivity radius is at least $1$. We denote by $\Lambda$ the bundle
of complex valued  
differential forms on $V$, and we keep this notation for its 
pull-back to $T^*V$ and its restrictions to various subsets of $V$ 
and $T^*V$.  For $x\in V$, we denote by $O_x$ the geodesic ball with radius 
$1/4$, $H_x$ the Hilbert space $L^2(O_x,\Lambda)$ and we write $H$ for 
the continuous field of Hilbert spaces  $\cup_{x\in V} H_x$. 
 
With a model operator on $\RR^n$, for instance those given in theorem 
(19.2.12)  of  \cite{hormander1-III}, we define a continuous 
family $P=(P_x)_{x\in V} $ of pseudodifferential operators 
$P_x\in\Psi^0(O_x,\Lambda)$ of order $0$ satisfying the following conditions~: 
\begin{enumerate} 
\item $P_x$ is {\sl trivial at infinity of } $O_x$, which means that 
$P_x$ is the sum of a compactly supported pseudodifferential operator 
and a smooth bounded section of the bundle $\End \Lambda \To 
O_x$, (in particular, this ensures boundedness on $H_x$), 
\item $P_x$ is selfadjoint on $H_x$, and has degree 
  one (i.e. $P_x=\begin{pmatrix} 0 & P_x^- \\ P_x^+ & 0 \end{pmatrix}$) 
  with respect to the grading  induced by $\Lambda=\Lambda^{ev}\oplus\Lambda^{odd}$, 
\item $P^2_x-Id$ is a compactly supported pseudodifferential operator 
of order $-1$; in particular it is compact on $H_x$, 
\item the family $P=(P_x)_{x\in V}$  has a 
trivial index bundle of rank one. In fact, for all $x$, $P_x^+$ is
onto, it has a one dimensional kernel and there exists a continuous
section $V\ni x\mapsto e_x\in \ker P_x^+\subset L^2(M,\Lambda)$. 
\end{enumerate} 
Here the continuity of the family means that $P$ is an endomorphism of the  $C(V)$-Hilbert
module $H$.  \\ We let $a_x$ be 
the principal symbol of $P_x$. Under the assumptions above, the
Kasparov module $$ \lambda_x = [(C_0(T^*O_{x},\Lambda),1,a_{x})]$$ 
is a generator of $K_0(C_0(T^*O_{x}))\simeq\ZZ$. The following element
: 
 $$ \lambda_V = [(C_0(T^*O_{x},\Lambda),1,a_{x})_{x\in V}]\in 
 K_0(C_0(T^*V\times V))$$ 
is  the dual Dirac element used in the proof the Poincar{\'e} duality 
between $C(V)$ and $C_0(T^*V)$ \cite{Ka2}. \\ There is an 
alternative elegant description of $\lambda_V$ \cite{CoS}. Let us consider 
the map $f : V\to T^*V\times V \ , \ x\mapsto ((x,0),x)$. This map 
is $K$-oriented so it gives rise to an element $f!\in KK(C(V), 
C_0(T^*V\times V))$. If $p$ 
denotes the obvious map $\CC \to C(V)$, then : 
 $$ \lambda_V = [p]\ot f! \ . $$ 

With the previous example in mind  and keeping the same notations,  we
shall define an element $\Delta\in K_0(C^{*}(\cG\times X))$ whose
evaluation $\lambda= (e_{0})_{*}(\Delta)\in K_0(C^{*}(G\times X))$
will be the appropriate dual Dirac element of  
the pseudomanifold $X$. 

We define a map $h :X\setminus\{c\} \simeq M^+ \longrightarrow M$ which pushes  
points in $M^{-}$. More precisely, $h(y)=y$ when $k_y\geq 1$,  
and $h(y)=( y_L,l(k_y))$ 
otherwise,  where : 
$$ l(k)=\begin{cases} 3k-2 & \text{ if } 1/2 \le k\le 1 \\ 
                     -1/2    & \text{ if } 0 < k\le  1/2  \ .\end{cases} 
 $$ \label{HM} 
>From now we fix $\varepsilon \in ]0, 1/2[$. Recall that
$X_\varepsilon=\{x\in X \ | \ k_x\geq \varepsilon \}$. We set : 
 $$ \delta = (\lambda_{h(x)})_{t\in [0,1],x\in X_\varepsilon}\in  
 K_0(C_0(\cA^*(\cG)\times X_\varepsilon)). $$  
 Here  $\cA^*(\cG)\simeq T^*M\times [0,1]$ and $ \delta$ corresponds
 to the $K$-oriented map :  
$[0,1]\times X_\varepsilon \to T^*M\times[0,1]\times X_\varepsilon,\
(t,x)\mapsto ((h(x),0),t,x)$.  
 
\noindent Next, let us consider the adiabatic 
groupoid of $\cG$ \cite{Co0,MP,NWX} (see example \ref{ex2} (1))~:  
$$\cH =\{0\}\times\cA(\cG) \cup ]0,1]\times \cG\rightrightarrows
[0,1]\times\cG^{(0)}.$$   
We let $\Theta \in KK(C_{0}(\cA^*(\cG)),C^{*}(\cG))$ be the  
$KK$-element associated to $\cH$, i.e.~: 
$$ \Theta = [ev_0]^{-1}\ot [ev_1], $$ 
where $ev_{0}: C^{*}(\cH)\rightarrow C_{0}(T^*M\times [0,1])$ is the evaluation 
map at $0$ composed with the Fourier transform 
$C^*(\cA(\cG))\overset{\simeq}{\to} C_0(\cA^*(\cG))$, and 
$ev_1 : C^*(\cH)\to C^*(\cG)$ is the evaluation at $1$. 
We define $\Delta_\varepsilon \in K_0(C^*(\cG\times X_\varepsilon))$ by : 
 $$ \Delta_\varepsilon= \delta \underset{C_0(\cA^*(\cG))}{\ot} \Theta\
 . $$ 
\begin{prop} 
\begin{enumerate} 
\item The element $\Delta_\varepsilon$ satisfies : 
    $$ (e_1)_*(\Delta_\varepsilon)=1_{X_\varepsilon}\in
    K^0(X_\varepsilon )\simeq K_0(\cK\otimes C(X_\varepsilon))$$
 where $e_1 : C^*(\cG)\to  C^*(\cG\vert_{M\times\{1\}})\simeq\cK$ is the evaluation map at $1$.
 \item There exists $\Delta_0\in K_0(C^*(\cG\times X |_{\cO\times[0,1]}))$ extending
  $\Delta_\varepsilon$, that is : 
 $$ r_*\circ (i_{\cO\times[0,1]})_*(\Delta_0)=\Delta_\varepsilon$$
where $\cO$ is the open subset $\cup_{x\in M^+}O_{h(x)}\times \{x\} \cup
M^-\times \co L $ of $M\times X$, $i_{\cO\times[0,1]}: C^*(\cG\times X |_{\cO\times[0,1]})
\longrightarrow C^*(\cG\times X)$ is the inclusion morphism and $r : C^*(\cG\times X)\longrightarrow
C^*(\cG\times X_{\varepsilon})$ is the 
restriction morphism. 
\end{enumerate}
\end{prop}
\begin{pf}
1) Let us note $\cO_\varepsilon = \cup_{x\in X_\varepsilon }O_{h(x)}\times 
\{x\}$. This is an open subset of $M\times X_\varepsilon$ the
unit space of the groupoid $G\times X_\varepsilon$. Let
$\widetilde{\cO_\varepsilon}$ be its lift to  $M\times[0,1]^2\times X_\varepsilon$ which is the unit space 
of the groupoid $\cH \times X_\varepsilon$.  We let $\Omega$ be the groupoid : 
  $$ \Omega=\cH\times X_\varepsilon |_{\widetilde{\cO_\varepsilon}} $$
In fact, $\Omega$ identifies with the adiabatic groupoid of $\Omega_{1}=\cG\times
X_\varepsilon |_{\cO_\varepsilon \times [0,1]}$. 

We shall use the pseudodifferential calculus on $\Omega$ to get an
explicit representant of $\Delta_\varepsilon$. 
The family $(a_{h(x)})_{x\in X_\varepsilon}$ depends smoothly on $x$ and defines
a symbol $a\in S^0(\cA^*(\Omega),\End\Lambda)$. Note that this symbol is independent of
the two real parameters coming from the lift of $\cO_\varepsilon$ to $\widetilde{\cO_\varepsilon}$.
Let  $\mbox{Op}_{\Omega}$ be a quantification map for $\Omega$.
Thanks to the properties of this calculus and the fact that 
each $a_{h(x)}(y,\xi)$ is of order $0$, trivial at infinity (that is
independent of $\xi$ near the infinity of $O_{h(x)}$) and
$a_{h(x)}^2(y,\xi)- 1$ is of order $-1$ and vanishes near the infinity
of $O_{h(x)}$, we  deduce from the exact sequence (\ref{E-S-PDO})~:  
 $$ \mbox{Op}_{\Omega}(a)\in \cL(C^*(\Omega,\Lambda))\text{ and } 
\mbox{Op}_{\Omega}^2(a)-\mbox{Id} \in  \cK(C^*(\Omega,\Lambda)).$$ 

Hence, we get an element 
$[(C^*(\Omega,\Lambda),\mbox{Op}_{\Omega}(a))]\in K_0(C^*(\Omega))$
which gives using the inclusion $\Omega\subset \cH\times X_\varepsilon$
an element $\widetilde{\Delta} \in K_0(C^*(\cH\times
X_\varepsilon))$ satisfying~: 
 $$ (ev_0)_*(\widetilde{\Delta})=\delta \in K_0(C_0(T^*M\times[0,1]\times X_\varepsilon))$$
Hence : 
 $$ \Delta_\varepsilon= (ev_1)_*(\widetilde{\Delta}) =
 [(C^*(\Omega_1,\Lambda),\mbox{Op}_{\Omega_1}(a))]$$ 
where  $\Omega_1=\cG\times X_\varepsilon |_{\cO_\varepsilon\times[0,1]}$. 
Now consider the evaluation map $e_1 :C^*(\cG)\to \cK $. We
get : 
 $$ (e_1)_*(\Delta_\varepsilon) =
 [(C^*(\Omega_{1,1},\Lambda),\mbox{Op}_{\Omega_{1,1}}(a))]\in
 K_0(\cK\otimes C(X_\varepsilon))$$
where we have set
$\Omega_{1,1}=\cG\times X_\varepsilon \vert_{\cO_{\varepsilon} \times \{1\}}$. Note that :
 $$\Omega_{1,1}=\bigcup_{x\in X_\varepsilon }O_{h(x)}\times
 O_{h(x)}\times\{x\}\subset (M\times M)\times X_\varepsilon $$
and $\mbox{Op}_{\Omega_{1,1}}$ is an ordinary quantification map which
assigns to a symbol living on $T^*O_{h(x)}\simeq \cA^*(O_{h(x)}\times
O_{h(x)})$ a pseudodifferential operator on $O_{h(x)}$. Since 
$P\vert_{X_\varepsilon}=(P_{h(x)})_{x\in X_\varepsilon}$ has symbol equal to
$(a_{h(x)})_{x\in X_\varepsilon}$ and has a trivial index bundle of rank one, the
following holds : 
$$ (e_1)_*(\Delta_\varepsilon) =
 [(C^*(\Omega_{1,1},\Lambda),\mbox{Op}_{\Omega_{1,1}}(a))]=
 [(H\vert_{X_\varepsilon},P\vert_{X_\varepsilon})]=1_{X_\varepsilon}\in K^0(X_\varepsilon) \ .$$

\noindent 2) The existence of $\Delta_0$ follows immediately from : 
\begin{lem} If $r_L : C(X_\varepsilon)\longrightarrow  C(L)$
  ($L=\partial X_\varepsilon $) denotes the restriction homomorphism, 
  and 
$i_{M^-\times[0,1]} : \cK \ot C([0,1])\simeq
C^*(\cG\vert_{M^-\times[0,1]}) \longrightarrow C^*(\cG)$
  the inclusion morphism,   
  then 
 $$(r_L)_*(\Delta_\varepsilon) = (i_{M^-\times[0,1]})_*(1_{L}) $$ 
 where $1_{L}$ is the unit of the ring $K_0(\cK\ot
 C([0,1]\times L))\simeq K^0(L)$.  
\end{lem} 
\begin{pf} 
The element $(r_L)_*(\Delta_\varepsilon)$ is represented by 
   $$ (C^*(\partial \Omega_1,\Lambda), \partial P)\in
   E(\CC,C^*(\cG\times L))$$ 
where $\partial \Omega_1 = \cup_{(t,x)\in [0,1]\times L}
O_{h(x)}\times O_{h(x)}\times \{(t,x)\}\subset (M^-\times M^-)\times[0,1]\times L$
and $\partial P=(P_{h(x)})_{x\in L}$. Since $C^*(\partial
\Omega_1,\Lambda)$ is also a $\cK\ot C([0,1]\times L)$-Hilbert module, we observe that : 
   $$ x= [(C^*(\partial \Omega_1,\Lambda), \partial P)]\in
   K_0(\cK\otimes C([0,1]\times L))$$
is such that $(i_{M^-\times[0,1]})_*(x)= (r_L)_*(\Delta_\varepsilon)$. Moreover, 
under the isomorphism $K^0(L)\simeq K_0(\cK\ot
 C([0,1]\times L))$, the element $x$ is represented by : 
$(H_{h(x)}, P_{h(x)})_{x\in L}$, which also represents  the unit
 element $1_{L}\in K^0(L)$ thanks to the triviality of the index bundle of
 the family $ (P_{h(x)})_{x\in L}$.
\hfill$\Box$\end{pf} 
By a slight abuse of notation, $\cG \times X_\varepsilon |_{\cO\times
  [0,1]}$, $\cG \times c_\varepsilon L |_{\cO\times [0,1]}$ and 
$\cG \times L |_{\cO\times [0,1]}$ will denote respectively the
restrictions of $\cG \times X$ to $\cO\times [0,1] \cap M\times
X_\varepsilon \times [0,1]$, $\cO\times [0,1] \cap M\times
c_\varepsilon L \times [0,1]=M^-\times
c_\varepsilon L \times [0,1]$ and $\cO\times [0,1] \cap M\times
\partial X_\varepsilon \times [0,1]=M^-\times L \times [0,1]$. 
 
It is obvious from the concrete
description of $\Delta_\varepsilon$ that it comes from an element 
$\Delta_{\varepsilon,\cO} \in K_0(C^*(\cG \times X_\varepsilon
|_{\cO\times [0,1]}))$ via the inclusion
morphism. Now let $x_{0}\in K_0(C^{*}(\cG \times
c_\varepsilon L |_{\cO\times [0,1]}))$  be the pushforward of
$1\in K_0(\CC)$ via the obvious homomorphism :  
$$ \cK\simeq C^{*}(M^{-}\times M^{-})\To C^{*}(M^{-}\times M^{-})  
\ot C([0,1] \times c_\varepsilon L) \simeq C^{*}(\cG \times c_\varepsilon L |_{\cO\times
  [0,1]})) .$$ 
The preceding lemma and the Mayer-Vietoris exact sequence in
$K$-theory associated to the following commutative diagram ~:   
$$
\begin{CD} 
C^*(\cG \times X |_{\cO\times [0,1]}  )       @>r>>  C^* (\cG \times
X_{\varepsilon } |_{\cO\times [0,1]} ) \\ 
 @Vr_{cL}VV             @VVr_{L}V \\ 
 C^* (\cG\times c_\varepsilon L |_{\cO\times [0,1]})   @>>r_{L}>
 C^* (\cG\times L |_{\cO\times [0,1]})\simeq \cK\ot C([0,1]\times L)
 \end{CD}
 $$ 
show that there exists $\Delta_{0}\in K_0(C^{*}(\cG\times X|_{\cO\times [0,1]}))$
satisfying $r_*(\Delta_0)=\Delta_{\varepsilon,\cO}$ and $(r_{cL})_*(\Delta_0)=x_0$.
\hfill$\Box$\end{pf} 
\begin{rem} In fact, $(r_L)_{*}(\Delta_{\varepsilon,\cO})$ may be represented by the {\sl trivial} vector bundle 
$\ker P\vert_{L}=\cup_{x\in L}\ker P_{h(x)}\To L$ while  $x_{0}$ may be represented by the {\sl product} vector bundle $cL\times\CC\To cL $. The element $\Delta_{0}$ is obtained by gluing these  bundles along $L=L\times\{1\}\subset cL$. This involves a bundle isomorphism 
$\ker P \vert_{L}\simeq L\times \CC$, which yields a continuous map 
$\psi : L\to GL_{1}(\CC)$ and a class  $[\psi]\in K^{1}(L)$. One could
be more careful with the construction of  the family $(P_x)$ to make
sure that $[\psi]=0$, otherwise one may perturb $\Delta_{0}$ by
elements coming from  $K^1(L)$. That will be done in
the next proposition. 
 \end{rem}

\begin{prop} There exists $\Delta_\cO\in K_0(C^*(\cG\times X |_{\cO\times[0,1]}))$ such
  that :  \\
1) $r_*\circ (i_{\cO\times [0,1]})_*(\Delta_\cO)=\Delta_\varepsilon $, \\
2) $(e_{1})_{*}\circ (i_{\cO\times [0,1]})_*(\Delta_\cO)=1_X\in K^0(X)\simeq K_0(C(X)\otimes\cK)$.
\end{prop} 
\begin{pf}
Firstly, we note that  : 
 $$ r_*\circ(e_{1})_{*}\circ (i_{\cO\times [0,1]})_*(\Delta_{0})  = (e_{1})_{*}\circ r_{*}\circ (i_{\cO\times [0,1]})_*(\Delta_0)   
 =  (e_{1})_{*} (\Delta_\varepsilon) = 1_{X_\varepsilon} =  r_*(1_X)  $$ 
>From the exact sequence : 
$$
\begin{CD} 
0 @>>> C_0(\co_\varepsilon L) @>j>> C(X)  @>r>> C(X_\varepsilon) @>>> 0
 \end{CD}
 $$ 
and the previous computation, we deduce : 
  $$(e_{1})_{*}(\Delta_0) - 1_X \in \mbox{Im}(j_*)$$ 
We choose $y_0\in K^1(L)\simeq K_0(C(\co_\varepsilon L)\otimes\cK)$ such that 
$$ j_*(y_0)=(e_{1})_{*}(\Delta_0) - 1_X$$ 
Moreover one deduces from the fact that the inclusion 
$\cK(L^2(M^-))\subset \cK(L^2(M))$ induces an isomorphism in
$K$-theory, that  
$$(e_1)_* \circ i_*: 
 K_0(C^*(\cG\times \co_\varepsilon L |_{\cO\times [0,1]}))\To  K_0(\cK\otimes C_0(\co_\varepsilon L))$$ 
is an isomorphism. Where $\cG\times \co_\varepsilon L
 |_{\cO\times [0,1]}$ is the restriction of $\cG \times X$ to
   $\cO\times [0,1] \cap M\times [0,1] \times \co_\varepsilon L$ and $i:
   C^*(\cG\times \co_\varepsilon L |_{\cO\times [0,1]}) \To C^*(\cG  
     \times X)$ is the inclusion morphism. 
 
Now let $\widetilde{y_0}\in K_0(C^*(\cG \times \co_\varepsilon L |_{\cO\times
  [0,1]}))$ be the unique element such that $(e_1)_* \circ
i_*(\widetilde{y_0})=-y_0$. We set~:  
 $$ \Delta_\cO = \Delta_0+j_*(\widetilde{y_0})  $$
where we have still denoted by $j$ the inclusion morphism from 
$C^*(\cG\times \co_\varepsilon L |_{\cO\times [0,1]})$ to $C^*(\cG\times X 
|_{\cO\times [0,1]})$. Then $\Delta_\cO\in K_0(C^*(\cG\times
X|_{\cO\times [0,1]}))$ satisfies (1) and (2).
\hfill$\Box$\end{pf}
We define $$\Delta=(i_{\cO\times [0,1]})_*(\Delta_\cO) \in
K_0(C^*(\cG \times X)) \ .$$
\begin{defn} The dual Dirac element $\lambda\in K_0(C^*(G) \ot C(X))$ of 
the singular manifold $X$ is defined by : 
    $$ \lambda = (e_0)_*(\Delta)$$
where $e_0 : C^*(\cG)\to C^*(G)$ is the evaluation homomorphism at $0$. 
\end{defn} 
We have proved the following : 
\begin{prop}\label{propproduit0} (1) The following equality holds : 
   $$ \lambda\underset{C^*(G)}{\otimes}\partial = 1_{X}\in 
   K_0(C(X))\simeq K^0(X),$$ 
where $1_{X}$ is the unit of the ring $K_0(C(X))$. \\
(2) For each  open subset $\cO_{\alpha}$,  $0<\alpha<1$, of $M\times X$ defined by : 
$$ \cO_{\alpha}= \{ (x,y)\in M\times X \ \vert \ d(x,y)<1, \ k_{y}>0 \} \cup 
\{ (x,y)\in M\times X \ \vert \ k_{x}< \alpha ,\ k_{y}< \alpha \}$$
the dual Dirac  element $\lambda$ belongs to the range of : 
$$(i_{\cO_{\alpha}})_{*}\ :\ K_0(C^{*}(G\times X\vert_{\cO_{\alpha}}))\To K_0(C^{*}(G\times X))$$ 
\end{prop} 
Roughly speaking,  every $\cO_{\alpha}$ contains the "support" of
$\lambda$. From now on, we choose $\cO=\cO_{1/2}$ and $\widetilde{\lambda}$ a preimage of $\lambda$ for $(i_{\cO})_{*}$.  
\section{The Poincar{\'e} duality} 
This section is devoted to the proof of our main result : 
\begin{thm}\label{thmedp} The Dirac element $D$ and the dual Dirac element $\lambda$ 
  induce a Poincar{\'e} duality between $C^*(G)$ and $C(X)$, that is~: 
\begin{enumerate} \item 
$ \lambda\underset{C^*(G)}{\otimes} D =  1_{C(X)}\in KK(C(X),C(X)).$ 
\item $ \lambda\underset{C(X)}{\otimes} D = 1_{C^*(G)}\in KK(C^*(G),C^*(G)).$ 
\end{enumerate} 
\end{thm} 
                         
%
\subsection{Computation of $\lambda {\otimes}_{C^*(G)} D$} 
\label{produit1} Let $m : C(X)\ot C(X) \rightarrow C(X)$ be 
the morphism of $C^*$-algebras induced by multiplication of 
functions. 
\begin{lem} The following equality holds : $$ 
  \lambda \underset{C^*(G)}{\ot} [\Psi] = \lambda 
  \underset{C(X)}{\ot} [m] \ .$$ 
\end{lem} 
\begin{pf} According to proposition \ref{propproduit0} we have that 
$\lambda \underset{C^*(G)}{\ot} [\Psi] = 
\tau_{C(X)}(\tilde{\lambda}) \ot [H_0]$ and 
$\lambda\underset{C(X)}{\ot} [m]= 
\tau_{C(X)}(\tilde{\lambda}) \ot [\tilde H_1]$, where 
  $H_0$ and $\tilde H_1$ are the morphisms from $C^*(G\times X \vert _{\cO})\ot C(X)$ 
to $C^*(G\times X)$ defined by : 
  $$H_0(B\ot f)(\gamma,y)= f(r(\gamma))B(\gamma,y)  
\ \mbox{ and } 
    \tilde H_1(B\ot f)(\gamma,y)=f(y)B(\gamma,y) \ ,$$ when 
$f \in C(X)$, $B\in C_c(G\times X \vert _{\cO})$, $(\gamma,y)\in 
G\times X$, and $C(X)$ is identified  with the algebra of 
continuous functions on $M$ which are constant on $M^-$. \\ There 
is an obvious homotopy between $\tilde H_1$ and $H_1$ defined by 
$$H_1(B\ot f)(\gamma,y)= f(h(y))B(\gamma,y) ,$$ where $h$ is the 
map constructed in section \ref{HM} to push points in $M^-$.\\ Let 
$\cW = \{ (x,y)\in M\times X \ ; \ d(x,y)< 1 \text{ and } k_y \ge 
\varepsilon \} $ viewed as a subset of $ M^+ \times M^+$. Let $c$ be 
the continuous function $\cW \times [0,1]\rightarrow M^+ \subset 
X$ such that $c(x ,y ,\cdot)$ is the geodesic path going from $x$ 
to $y$ when $(x,y)\in \cW$. \\ We obtain an homotopy between $H_0$ 
and $H_1$ by setting for each $t\in [0,1]$ : $$H_t(B\ot 
f)(\gamma,y)=\left\{ 
  \begin{array}{ll} f(c(r(\gamma),h(x),t))B(\gamma,x) \mbox{ if } k_x 
    \geq 1/2 \\ f(c)B(\gamma,x) \mbox{ elsewhere} . \end{array} \right. $$ 
\hfill $\Box$ \end{pf} 
Now we are able to compute the product $\lambda \underset{C^*(G)}{\ot}D $ : 
$$\begin{array}{ccccc} \lambda \underset{C^*(G)}{\ot} D 
& = &  (\lambda \underset{C^*(G)}{\ot} [\Psi]) \underset{C^*(G)}{\ot} 
\partial   &  = & (\lambda \underset{C(X)}{\ot} [m])\underset{C^*(G)}{\ot} 
\partial \\  & = & 
\tau_{C(X)}(\lambda) \ot ([m] \underset{\CC}{\ot} 
\partial)   & = & \tau_{C(X)}(\lambda) \ot (\partial \underset{\CC}{\ot} 
[m]) \\ 
& = & \tau_{C(X)}(\lambda 
\underset{C^*(G)}{\ot} \partial ) \ot [m] & =  & \hspace{-2cm} 1_{C(X)}\ . 
\end{array}$$ 
\noindent The equality of the first line results from the previous 
lemma, the equality of the second line comes from the 
commutativity of the Kasparov product over $\CC$ and the last 
equality follows from proposition \ref{propproduit0}. This 
finishes the proof of the first part of theorem \ref{thmedp}. 
 
\smallskip \noindent Let us notice the following consequence : for every $C^*$-algebras $A$ and $B$, we 
have :
$$(\cdot \underset{C^*(G)}{\ot} D) \circ (\lambda 
  \underset{C(X)}{\ot} \cdot) = Id_{KK(C(X)\ot A,B)} \ \mbox{ and} \label{iso1}\\ 
    (\lambda \underset{C^*(G)}{\ot} \cdot) \circ (\cdot \underset{C(X)}{\ot} 
  D)=Id_{KK(A,B\ot C(X))} .$$ 
In particular, this implies the following useful remark : 
\begin{rem} \label{rem} The equality 
$(\cdot\underset{C^*(G)}{\ot}D)(\lambda\underset{C(X)}{\ot}D)=(\cdot\underset{C^*(G)} 
{\ot}D)\circ(\lambda\underset{C(X)}{\ot}\cdot)(D)=D$ ensures that 
$\lambda \underset{C(X)}{\otimes} D-1_{C^*(G)}$ belongs to the 
kernel of the map $(\cdot \underset{C^*(G)}{\ot} D) : 
KK(C^*(G),C^*(G)) \rightarrow KK(C^*(G)\ot C(X),\CC) \ .$
\end{rem} 
\subsection{Computation of $\lambda {\otimes}_{C(X)}  D$} 
The purpose here is to prove that 
\begin{eqnarray*} \lambda \underset{C(X)}{\otimes} 
  D=\tau_{\underline{C^*(G)}}(\lambda)\ot \tau_{C^*(G)}([\Psi]\ot \partial)=1_{C^*(G)}. 
 \end{eqnarray*} 
This problem leads us to study the invariance of 
$\lambda\underset{C(X)}{\ot} [\Psi]$ under the {\it flip} 
automorphism $\tilde f$ of $C^*(G\times G)\simeq C^{*}(G)\ot C^{*}(G)$.
 
\noindent  Indeed, we have : 
$$  \tau_{C^*(G)}([\Psi]) \ot [\tilde f]\ot  \tau_{C^*(G)}(\partial)= [\Psi] \underset{\CC}{\ot} 
\partial =\partial \underset{\CC}{\ot}[\Psi]=\tau_{\underline{C^{*}(G)}\otimes C(X)}(\partial) \ot [\Psi], $$ 
which implies (cf. proposition \ref{propproduit0}) : 
\begin{equation*} 
((\lambda\underset{C(X)}{\otimes}[\Psi]){\ot}[\tilde{f}]){\otimes}\tau_{C^{*}(G)}(\partial) 
=  (\lambda  \underset{C^*(G)}{\ot}  \partial)\underset{C(X)}{\ot}[\Psi] =  1_{C^*(G)}. 
\end{equation*} 
Hence \begin{equation}  \label{invflip} \lambda 
\underset{C(X)}{\otimes}  D - 1_{C^*(G)}= 
((\lambda\underset{C(X)}{\otimes}[\Psi]){\ot}([\id]-[\tilde{f}])){\otimes}\tau_{C^{*}(G)}(\partial) 
\end{equation}  and the invariance of 
$\lambda\underset{C(X)}{\otimes}[\Psi]$ under $[\tilde{f}]$  would 
enable us to conclude the proof. Such an invariance would be 
analogous to lemma 4.6 of \cite{Ka2}. Unfortunately we are not 
able to prove that $\lambda\underset{C(X)}{\ot} [\Psi]$ 
is invariant under the flip. 
 
\medskip \noindent 
Put $C=L\times ]-1,1[$ and $F=M\times M \setminus C\times C$.\\ We 
let $\tilde{\cO}$ be the inverse image of $\cO$ by the canonical 
projection of $M\times M \rightarrow M\times X$. We denote 
$(i_{\tilde{\cO}})_{*} : KK(C^*(G), C^*(G\times G 
\vert_{\tilde{\cO}}))  \rightarrow KK(C^*(G),C^*(G\times G))$ the 
morphism corresponding to the inclusion $i_{\tilde{\cO}}$. A 
simple computation shows that : 
 \begin{lem}\label{lem2} The element $\lambda \underset{C(X)}{\ot} [\Psi] $ belongs 
  to the image of $(i_{\tilde{\cO}})_*$. 
\end{lem} 

 The set $F \cap \tilde{\cO}$ is an open and symmetric subset of $M\times M$, hence the flip makes sense on 
 $C^{*}(G\times G \vert _{F \cap\tilde{\cO}})$. 
\begin{lem} 
The flip automorphism of $C^*(G\times G \vert_{ F \cap \tilde{\cO}})$ 
is homotopic to identity. 
\end{lem} 
\begin{pf} The set $F\cap\tilde{\cO}=\{(x,y)\in M\times M\ \vert\ d(x,y)<1\ ,\  k_x\geq1\ \mbox{ or }\ k_y\geq 1\}$ 
is a subset of $\Mp \times \Mp$, so the algebra 
$C^*(G\times G \vert_{F \cap \tilde{\cO} })$ is isomorphic to $C_0(T^*(F \cap \tilde{\cO} ))$. 
To prove the lemma, it is sufficient to find a proper homotopy between the flip 
  $f_{ F \cap \tilde{\cO}} $ of $T^{*}(F \cap \tilde{\cO} )$ and $id_ {T^{*}(F \cap \tilde{\cO} )}$.\\ 
  The exponential map 
  of $M$ provides an isomorphism $\phi$ between $T^*(F \cap \tilde{\cO} )$ and 
  $(T^*M_{l})^{\oplus 3}$, where $M_{l}=\{x\in M \ \vert \ k_{x}\geq l\}$ for some $0<l<1$. 
  Via this isomorphism, the flip becomes the automorphism of 
  $C_0((T^*M_{l})^{\oplus 3})$ induced by 
  the map $g:(x,X,Y,Z) \in (T_x^*M_{l})^3 \mapsto (x,-X,Z,Y)$. One can take 
  for example : 
\centerline{ $\phi:(x,y,X,Y) \mapsto (m(x,y), 
exp_{m(x,y)}^{-1}(x)- exp_{m(x,y)}^{-1}(z),T(x,y,X),T(y,x,Y)),$} 
where 
  $m(x,y)=exp_x(\frac{exp_x^{-1}y}{2})$ is the middle point of the geodesic 
  joining $x$ to $y$ and $T(x,y,\cdot) : T_x M \rightarrow T_{m(x,y)} 
  M$ is the parallel transport along the geodesic joining $x$ to 
  $m(x,y)$.\\ Let $A:[0,1] \rightarrow SO_3(\RR)$ be a 
  continuous path from $\begin{pmatrix} -1 & 0 & 0 \\ 0 & 0 & 1 \\ 0 & 
    1 & 0 \end{pmatrix}$ to $\begin{pmatrix} 1 & 0 & 0 \\ 0 & 1 & 0 \\ 0 & 
    0 & 1 \end{pmatrix}$.\\ 
The map $[0,1]\times (T^*M_{l})^{\oplus 3} \rightarrow 
(T^*M_{l})^{\oplus 3}$; $(t,x,V)\mapsto (x, A_t . V)$ 
is a proper homotopy between identity and $g$. \hfill $\Box$ 
\end{pf} 
Note that $C\times C$ is a saturated open subset of $(G\times 
G)^{(0)}$. So we obtain the following commutative diagram of 
$C^*$-algebras : 
\begin{equation} \label{eqiCC} 
\begin{CD} 0\rightarrow @. C^*(G\times G \vert_{ C\times C \cap \tilde{\cO}}) 
  @>>>C^*(G\times G \vert_{\tilde{\cO}}) @>r_{F\cap \tilde{\cO}}>> 
  C^*(G\times G \vert_{ F \cap \tilde{\cO}}) @. \rightarrow 0\\ &&@VVV 
  @VVi_{\tilde{\cO}}V @VVi_{F\cap \tilde{\cO}}V \\  0\rightarrow @. C^*(G\times G 
  \vert_{C\times C}) @>>i_{C\times C}>C^*(G\times G) @>>r_{F}> 
  C^*(G\times G \vert_{F}) @. \rightarrow 0 \end{CD}  \end{equation} 
 Since flip automorphisms commute with restriction and 
inclusion morphisms, this commutative diagram and the previous 
lemma imply that the induced morphisms of $KK$ groups satisfy : 
$$(r_{F})_* \circ \tilde{f}_* \circ (i_{\tilde \cO})_* = 
(r_{F})_*\circ (i_{\tilde \cO})_* \ .$$ In other words, $(r_F)_* 
\circ (id-\tilde{f})_* \circ (i_{\tilde \cO})_* $ is the zero map. 
\\ Hence lemma \ref{lem2}  implies that 
$(\lambda\underset{C(X)}{\ot}[\Psi])\ot({[\id]-[\tilde f]})$ 
belongs to  the kernel of the map $(r_F)_*:KK(C^*(G),C^*(G\times 
G))\rightarrow KK(C^*(G),C^*(G\times G\vert_F))$. It follows from 
the long exact sequence in $KK$-theory associated to the second 
short exact sequence of (\ref{eqiCC}) 
that  $(\lambda\underset{C(X)}{\ot}[\Psi])\ot({[\id]-[\tilde f]})$ 
belongs to the image of the map 
$(i_{C\times C})_*:KK(C^*(G),C^*(G\times G\vert _{C\times C}))\rightarrow KK(C^*(G),C^*(G\times G))$. 
\begin{rem}\label{remkk} The $C^{*}$-algebra $C^*(G\vert _{C})$ is $KK$-equivalent to 
$\cK(L^2(M^-))$. Indeed, we have the following exact sequence~: 
$$ 
\begin{CD} 0 \to  \cK(L^2(M^-)) 
@>{i_{M^-,C}}>>  C^*(G \vert _{C}) @>>> C^*(G 
\vert _{C\setminus M^-}) \simeq 
C_0(T^*( L\times [0,1[)) \to 0 
\end{CD} 
$$ and the $C^*$-algebra $C_0(T^*( L\times [0,1[))$ is 
contractible. So $[i_{M^-,C}]$ is an invertible element of 
$KK(\cK(L^2(M^-)),C^*(G \vert _{C}))$.\\ 
 In particular, $C^*(G\times G\vert _{C\times C})$ is $KK$-equivalent to $\cK \ot \cK$. Furthermore, 
 the flip automorphism of $\cK \ot \cK$ is homotopic to identity. Together with our last result, this only 
 shows that $(\lambda\underset{C(X)}{\ot}[\Psi])\ot({[\id]-[\tilde f]})$ is a torsion element (of order $2$). 
\end{rem} 
\begin{lem}\label{propp} The element $\lambda\underset{C(X)}{\otimes} D - 1_{C^*(G)}$ 
belongs to 
the image of the map $(i_{C})_*:KK(C^*(G),C^*(G\vert_C)) \rightarrow KK(C^*(G),C^*(G))$ 
induced by the inclusion morphism $i_C$.\end{lem} 
\begin{pf} The proof follows from the equality $[i_{C\times C}] \ot \tau_{C^*(G)}(\partial)= 
([i_{C}] \underset{\CC}{\ot}[i_{C}] )\ot \tau_{C^*(G)}(\partial) 
=\tau_{C^*(G\vert_{C})}([i_{C}] \ot \partial)\ot [i_{C}]$ and the 
fact that $(\lambda \underset{C(X)}{\ot}[\Psi])\ot([\mbox{Id}]- 
[\tilde{f}])$ is in the image of $(i_{C\times C})_{*}$. \hfill 
$\Box$ \end{pf} Thus, with the remark \ref{rem} in mind, it 
remains to show  that the map $(\cdot\underset{C^*(G)}{\ot} 
D)\circ (i_C)_* $ is injective. \\ We consider the morphisms 
$i^{-}: \cK(L^2(M^-)) \rightarrow C^*(G)$ and $i^\cK 
:\cK(L^2(M^-)) \rightarrow \cK(L^2(M))$ induced by the inclusion 
of functions. Since $i^\cK$ preserve the rank of operators, 
$(i^\cK)_*$ is an isomorphism. We let $e_c: C(X) \rightarrow \CC$ 
be the evaluation map at $c$. The map $e_c$ admits a right 
inverse, so $e_c^*$ is injective. 
\begin{prop}\label{diagprems} For any $C^*$-algebra $A$, the following 
  diagram is commutative 
$$\begin{CD} KK(A,\cK(L^2(M^-))) @>(i^-)_*>> KK(A,C^*(G)) \\ 
@V(\cdot\ot [i^\cK]\ot b)VV  @VV(\cdot \underset{C^*(G)}{\ot} D)V \\ 
     KK(A,\CC) @>>e_c^*> KK(A\ot 
    C(X),\CC) \end{CD}$$ 
\end{prop} 
\begin{pf} 
\smallskip \noindent For any  $x\in KK(A,\cK(L^2(M^-)))$ we write : 
$$(\cdot \underset{C^*(G)}{\ot} D)\circ 
  (i^-)_*(x)  =  \tau_{C(X)}(x\ot [i^-])\ot D   =  \tau_{C(X)}(x) \ot 
  \tau_{C(X)}([i^-]) \ot [\Psi]\ot \partial \ . $$ 
If $f\in C(X)$ and $k\in KK(L^2(M^-))$, we observe that 
$\Psi((i^-(k)\ot f))=f(c)i^-(k)=e_c(f)i^-(k)$. In particular 
$\tau_{C(X)}([i^-])\ot [\Psi]=\tau_\cK ([e_c])\ot [i^-]$. It 
follows that $$(\cdot \underset{C^*(G)}{\ot} D)\circ(i^-)_*(x)= 
\tau_{C(X)}(x) \ot \tau_\cK ([e_c]) \ot [i^-] \ot \partial . $$ 
Furthermore, the following commutative diagram of $C^*$-algebras : 
$$\begin{CD} \cK(L^2(M^-)) @<ev_1 \ot id<< 
  C^*(\cG \vert _{M^-\times [0,1]})\simeq C([0,1])\ot \cK(L^2(M^-)) 
@>ev_0 \ot id>>  \cK(L^2(M^-)) \\  @V{i^\cK}VV @VVV  @VVi^-V \\ 
\cK(L^2(M)) @<<e_1< C^*(\cG )  @>>e_0> C^*(G) \end{CD} $$ shows 
that $[i^-] \ot \partial= [i^\cK]$. \\ Finally, using that 
$\tau_{C(X)}(x) \ot \tau_\cK ([e_c])= x\underset{\CC}{\ot}[e_c]$, 
we get $$\begin{array}{ccc} (\cdot \underset{C^*(G)}{\ot} D)\circ 
  (i^-)_*(x) & = & e_c^*\circ (\cdot \ot [i^\cK]\ot b )(x). \end{array}$$ 
\hfill $\Box$ \end{pf} 
 We have already noticed that $C^*(G\vert _{C})$ is 
$KK$-equivalent to $\cK(L^2(M^-))$ (cf. remark \ref{remkk}), and 
$i_{C}\circ i_{M^-,C}= i^-$. So, using the previous 
proposition (applied to $A=C^*(G)$), we deduce :
\begin{cor}\label{dercor} The morphism $(\cdot \underset{C^*(G)}{\ot} D)$ is injective when 
restricted to the image of $(i_C)_*$ going from 
$KK(C^*(G),C^*(G\vert_C))$ to $KK(C^*(G),C^*(G))$. 
\end{cor} 
 Combining lemma \ref{propp}, remark \ref{rem} and 
corollary \ref{dercor}, we conclude that $$\lambda 
\underset{C(X)}{\ot} D= 1_{C^*(G)}.$$ This finishes the proof of 
theorem \ref{thmedp}. 
\begin{rem} The $K$-duality for the pseudomanifold $X$ is strongly
  related to a Poincar{\'e} duality for manifolds with boundary.  
Let us consider the following two exact sequences : 
\begin{eqnarray} \label{suite0} 
       \begin{CD} 0 @>>> \cK(L^2(M^-)) @>{i^-}>>  C^*(G) @>{r}>> 
                 C_0(T^*\overline{M^+}) @>>> 0 
       \end{CD} 
\end{eqnarray} 
\begin{eqnarray} \label{suite} 
     \begin{CD} 0 @<<< \CC @<<{e_c}< C(X) @<<{j}< C_0(X\setminus \{c\}) 
               @<<< 0 
     \end{CD}   
\end{eqnarray} 

Note that the exact sequence (\ref{suite}) is split, and that 
proposition \ref{diagprems} ensures  the injectivity of $(i^-)_*$. Hence 
both (\ref{suite0}) and (\ref{suite}) give rise to short exact 
sequences in $KK$-theory, and invoking again proposition 
\ref{diagprems} , we get the following commutative diagram : 
\begin{equation*}\label{diagrho} 
\begin{CD}  0 \to KK(A,\cK ) @>{(i^-)_*}>> 
  KK(A,C^*(G)) @>{r_*}>> KK(A,C_0(T^*\overline{M^+})) \to 0 \\ 
  @V{(\cdot \ot [i^\cK ]\ot b)}VV @VV(\cdot \underset{C^*(G)}{\ot} D)V 
  @VVV    \\ 
 0 \to 
  KK(A,\CC) @>>{e_c^*}>  KK(A\ot C(X),\CC) @>>{j^*}> KK(A\ot C_0(M^+),\CC) \to 0 
\end{CD} 
\end{equation*} 
The vertical arrows are isomorphisms : it is obvious for the left 
one and a consequence of theorem \ref{thmedp} for the middle one. 
Hence there is an induced isomorphism 
$KK(A,C_0(T^*\Mp))\to  KK(A\ot C_0(M^+),\CC)$ making the diagram commutative.

Conversely, using \cite{CoS}, one can prove the $K$-duality  
 between $C_0(T\overline M^+)$ and $C_0(M^+)$, and obtain from this an
 alternative proof of theorem \ref{thmedp}. This will be used in a
 forthcoming paper to extend this work to general pseudomanifolds. 
\end{rem} 

\bibliographystyle{plain} 
\bibliography{biblioconeJFA}

\end{document}